\documentclass[a4paper,10pt]{article}

\usepackage[ascii]{inputenc}
\usepackage{amsmath,amsfonts,amssymb,amsthm}

\newtheoremstyle{theorem}{5pt}{5pt}{\itshape}{}{\bfseries}{.}{.5em}{}

\theoremstyle{theorem}
\newtheorem{theorem}{Theorem}
\newtheorem{lemma}[theorem]{Lemma}
\newtheorem{corollary}[theorem]{Corollary}
\newtheorem{proposition}[theorem]{Proposition}

\usepackage{titlesec}
\titleformat{\section}
{\Large\bfseries}{\thesection}{1em}{\large}
\titlespacing*{\section}{0pt}{3.5ex plus 1ex minus .2ex}{2.3ex plus .2ex}

\titleformat{\subsection}
{\normalfont\bfseries}{\thesubsection}{1em}{\normalsize}
\titlespacing*{\section}{0pt}{3.5ex plus 1ex minus .2ex}{2.3ex plus .2ex}

\allowdisplaybreaks[3]

\begin{document}

\title{Moments and Oscillations of Exponential \\ Sums Related to Cusp Forms}

\author{Esa V\!. Vesalainen}

\date{}

\maketitle

\begin{abstract}
We consider large values of long linear exponential sums involving Fourier coefficients of holomorphic cusp forms. The sums we consider involve rational linear twists $e(nh/k)$ with sufficiently small denominators. We prove both pointwise upper bounds and bounds for the frequency of large values. In particular, the $k$-aspect is treated. As an application we obtain upper bounds for all the moments of the sums in question. We also give the asymptotics with the right main term for fourth moments.

We also consider the mean square of very short sums, proving that on average short linear sums with rational additive twists exhibit square root cancellation. This result is also proved in a slightly sharper form.

Finally, the consideration of moment estimates for both long and short exponential sums culminates in a result concerning the oscillation of the long linear sums. Essentially, this result says that for a positive proportion of time, such a sum stays in fairly long intervals, where its order of magnitude does not drop below the average order of magnitude and where its argument is in a given interval of length $3\pi/2+\varepsilon$ and so can not wind around the origin.
\end{abstract}

\section{Introduction}

\subsection{Linear exponential sums related to cusp forms}

Let $F$ be a fixed holomorphic cusp form of even weight $\kappa\in\mathbb Z_+$ for the full modular group $\mathrm{SL}(2,\mathbb Z)$. Then $F$ has the usual normalized Fourier expansion
\[F(z)=\sum_{n=1}^\infty a(n)\,n^{(\kappa-1)/2}\,e(nz),\]
where $z\in\mathbb C$ with $\Im z>0$.
It is interesting to study the linear exponential sums
\[A(M,\alpha)=\sum_{n\leqslant M}a(n)\,e(n\alpha),\]
and more generally
\[A(M,\Delta;\alpha)=\sum_{M\leqslant n\leqslant M+\Delta}a(n)\,e(n\alpha),\]
where $M\in\mathbb R_+$ is large, $\Delta\in\left[1,M\right]$, and $\alpha\in\mathbb R$.
These sums provide an interesting window into the life of the Fourier coefficients.
They are also connected to various other problems; as an example we mention second moments of the corresponding $L$-function (see e.g. the introduction to \cite{Miller2006}).

When $\Delta$ is small compared to $M$, the resulting short sums provide a natural analogue of the classical problems of analytic number theory studying various error terms in short intervals. Furthermore, good estimates for short sums also provide a practical tool for reducing smoothing error thereby potentially leading to better estimates (as examples we mention \cite{Jutila1987b, ErnvallHytonen2008} and the proof of Theorem \ref{pointwise-estimate} below). Pointwise bounds for short sums have been obtained by Jutila \cite{Jutila1987b}, and Ernvall-Hyt\"onen and Karppinen \cite{ErnvallHytonenEtKarppinen2008, ErnvallHytonen2008}.

The first estimate for long linear exponential sums involving holomorphic cusp form coefficients was proved by Wilton \cite{Wilton1929} in the course of proving an analogue of Voronoi's summation formula for cusp form coefficients (see Theorem \ref{voronoi-summation-formula} below) and that the $L$-function connected to the Ramanujan $\tau$-function has infinitely many zeros on the critical line. Wilton's estimate is that the long linear sum is $\ll M^{1/2}\,\log M$,
uniformly for $\alpha\in\mathbb R$.
Rankin \cite{Rankin1939} and Selberg \cite{Selberg1940} famously proved that
\[\int\limits_0^1\left|\sum_{n\leqslant M}a(n)\,e(n\alpha)\right|^2\mathrm d\alpha
=\sum_{n\leqslant M}\left|a(n)\right|^2=A\,M+O(M^{3/5}),\]
where $A$ is a certain positive real constant depending on $F$. This implies that only the logarithm might possibly be removed from Wilton's estimate, and this was done by Jutila \cite{Jutila1987b}. Indeed, we have
\[\sum_{n\leqslant M}a(n)\,e(n\alpha)\ll M^{1/2},\]
uniformly for $\alpha\in\mathbb R$.

\subsection{Linear exponential sums with rational twists}

In the study of linear exponential sums, the case in which $\alpha$ is near a fraction with a small denominator is often different in character from the one in which $\alpha$ is not close to such a fraction. Furthermore, the behavior near such rational values is often strongly linked to the behavior at such rational points, and it is this latter behavior that we are concerned with in this paper.

For rational values of $\alpha$ the study of the linear exponential sums is more akin to the classical problems of understanding the error terms in the Dirichlet divisor problem and the circle problem. Indeed, the analogies go rather deep (see e.g. \cite{Jutila1987a}), largely due to the fact that the divisor and circle problems, too, have modular origins.

However, the cusp form problems are in some respects more challenging. For example, the sums related to cusp forms have not been directly connected to such theories as the machinery of exponent pairs. Many of the best estimates are weaker than the best results on the analogous classical problems. An example of this is provided by the case $\alpha=0$: for the sums of Fourier coefficients, the best estimate to date, due to Rankin \cite{Rankin1989}, is
$A(M,0)\ll M^{1/3}\,\log^{-\delta}M$, for a certain small positive $\delta\in\mathbb R_+$. For comparison, during the 20th century the error term in the Dirichlet divisor problem saw a long string of upper bounds improving the cubic root cancellation first obtained by Voronoi \cite{Voronoi1903}. For more on the rather extensive literature on the Dirichlet divisor problem, we recommend Tsang's survey \cite{Tsang2010} and references therein, and Chapter 13 in the book \cite{Ivic2003} of Ivi\'c.

\subsection{What we do in this paper}

On average the long linear sum with rational twist $\alpha=h/k$ has the size $k^{1/2}\,M^{1/4}$, which is Theorem 1.2 in \cite{Jutila1987a}, and recalled in Theorem \ref{mean-square-of-long-sums} below. However, the pointwise upper bounds are far from this. The upper bound $k^{2/3}\,M^{1/3+\varepsilon}$ is Corollary on p.\ 30 in \cite{Jutila1987a}, and it is difficult to improve upon this, but we give a slight improvement for $M^{1/10}\ll k\ll M^{3/8-\varepsilon}$ by arguing from the Voronoi-type summation formula along the lines of Ivi\'c's paper \cite{Ivic2004} and reducing the smoothing error using the estimates for short exponential sums by Ernvall-Hyt\"onen and Karppinen \cite{ErnvallHytonenEtKarppinen2008}.

Using arguments of Ivi\'c's paper \cite{Ivic1983}, we prove in Theorem \ref{large-values} that, in a certain sense, a long linear exponential sum can not exceed its average value $k^{1/2}\,M^{1/4}$ too often.
These estimates will lead, as in \cite{Ivic1983}, to rather general moment estimates for such sums. In particular, the sums of coefficients will exhibit in every moment cancellation beyond the cubic root cancellation.

We also take a closer look at the fourth moment of those sums, obtaining the main term of the asymptotics, following Ivi\'c and Sargos \cite{IvicEtSargos2007}. The main term for the fourth moment for the divisor problem was first obtained by Tsang \cite{Tsang1992}.

We shall also consider the mean square of short sums with rational additive twists, proving essentially the expected square root cancellation on average:
\[\int\limits_M^{2M}\left|\sum_{x\leqslant n\leqslant x+\Delta}a(n)\,e\!\left(\frac{nh}k\right)\right|^2\mathrm dx\ll M\,\Delta\]
for $k\ll\Delta^{1/2-\varepsilon}$,
and
\[\int\limits_M^{2M}\max_{0\leqslant U\leqslant\Delta}\left|\sum_{x\leqslant n\leqslant x+U}a(n)\,e\!\left(\frac{nh}k\right)\right|^2\mathrm dx\ll M\,\Delta\,\log^2M\]
for $k\ll\Delta^{1/4-\varepsilon}$.
For sums of coefficients, the first of these essentially follows from the arguments of Jutila \cite{Jutila1984}, which we will follow. The mean square of square root length sums has been considered by Ernvall-Hyt\"onen \cite{ErnvallHytonen2011} who obtained square root cancellation with rational additive twists. The second estimate, which has already been obtained without twists by Wu and Zhai \cite{WuZhai2013}, is proved following Heath-Brown and Tsang as in the proof of Lemma 2 of \cite{Heath-BrownEtTsang1994}, and following the proof of Theorem 2 in \cite{IvicEtZhai2012}.
We would also like to mention that recently Ernvall-Hyt\"onen \cite{ErnvallHytonen2013} has considered the mean square of short exponential sums for which $\Delta$ is larger than $M^{1/2}$.

Our final Theorem \ref{oscillation-result} will be an analogue of the main theorem of Heath-Brown and Tsang in \cite{Heath-BrownEtTsang1994}, and it is the original motivation for the various moment estimates of this paper. Combining the various moment estimates for both long and short sums, we will see that, for $k\ll M^{1/8-\varepsilon}$, as $M$ increases, the long linear sum with a rational twist spends a positive proportion of its time in nearly square root length intervals in which its absolute value is at least of the average size $k^{1/2}\,M^{1/4}$ and where its argument is in a given fixed interval of length $3\pi/2$. Thus, essentially the result says that for a positive proportion of time, such sums do not exhibit much oscillation: the absolute value stays at least at the average level and the value does not wind around the origin.

\subsection{Notation}

The cusp form $F(z)$ will be a fixed holomorphic cusp form of even weight $\kappa\in\mathbb Z_+$ for the full modular group $\mathrm{SL}(2,\mathbb Z)$, and $F$ is assumed not to vanish identically. Its Fourier expansion will be normalized as follows:
\[F(z)=\sum_{n=1}^\infty a(n)\,n^{(\kappa-1)/2}\,e(nz).\]

The symbols $\ll$, $\gg$, $\asymp$ and $O$ are used for the usual asymptotic notation: for complex valued functions $f$ and $g$ in some set $\Omega$, the notation $f\ll g$ means that $\left|f(x)\right|\leqslant C\left|g(x)\right|$ for all $x\in\Omega$ for some implicit constant $C\in\mathbb R_+$. When the implicit constant depends on some parameters $\alpha,\beta,\ldots$, we use $\ll_{\alpha,\beta,\ldots}$ instead of mere $\ll$. The notation $g\gg f$ means $f\ll g$, and $f\asymp g$ means $f\ll g\ll f$.
The implicit constants are allowed to depend on the cusp form $F$ under study and $\varepsilon$.
The symbol $\varepsilon$ denotes an arbitrarily small positive real constant, and its value can and will change from one instance to the next. We emphasize that in the exponents of assumptions and conclusions the values of $\varepsilon$ can be chosen to be arbitrarily small.

Finally, if $f$ is a real-valued function, then
\[f_+=\max\{f,0\}.\]
If $f$ is complex-valued, then we will write
\[f_{++}=(\Re f)_++i(\Im f)_+.\]

\section{The results}

The classical upper bound \cite[Cor., p. 30]{Jutila1987a} which follows from the truncated Voronoi identity by taking absolute values is
\[\sum_{n\leqslant M}a(n)\,e\!\left(\frac{nh}k\right)\ll k^{2/3}\,M^{1/3+\varepsilon},\]
where $h$ and $k$ are coprime integers with $1\leqslant k\ll M$.
On the other hand, the estimates from \cite{Jutila1987b} give the upper bound $\ll M^{1/2}$, for arbitrary $k$.
It turns out that in the range $M^{1/10}\ll k\ll M^{3/8-\varepsilon}$, these two estimates can be improved upon by combining the argument of \cite{Ivic2004} with the short sum estimates from \cite{ErnvallHytonenEtKarppinen2008}.
\begin{theorem}\label{pointwise-estimate}
Let $M\in\left[1,\infty\right[$, and let $h$ and $k$ be coprime integers with $1\leqslant k\ll M^{1/2}$.
We have
\[\sum_{n\leqslant M}a(n)\,e\!\left(\frac{nh}k\right)
\ll\left\{\!\!\begin{array}{ll}
k^{2/3}\,M^{1/3}\,\log^{1/3}M&\text{when $k\ll M^{1/2}$,}\\
k^{1/4}\,M^{3/8+\varepsilon}&\text{when $M^{1/10}\ll k\ll M^{1/4}$,}\\
k^{2/3}\,M^{13/48+\varepsilon}&\text{when $M^{1/4}\ll k\ll M^{19/64}$,}\\
M^{15/32+\varepsilon}&\text{when $M^{19/64}\ll k\ll M^{21/64}$,}\\
k^{2/3}\,M^{1/4+\varepsilon}&\text{when $M^{21/64}\ll k\ll M^{3/8}$, and}\\
M^{1/2}&\text{when $k\ll M^{1/2}$.}
\end{array}\right.\]
\end{theorem}

The following theorem estimates the rarity of large values a long linear sum can take.
\begin{theorem}\label{large-values}
Let $M\in\left[1,\infty\right[$, and let $h$ and $k$ be coprime integers with $1\leqslant k\ll M^{1/2-\varepsilon}$, and let $V\in\mathbb R_+$ with $k\,M^\varepsilon\ll V\ll\sqrt M$.
Consider points $M\leqslant t_1<t_2<\ldots<t_R\leqslant2M$ satisfying $\left|t_r-t_s\right|\geqslant V$ for $r,s\in\left\{1,2,\ldots,R\right\}$ with $r\neq s$. Fix an exponent pair $\left\langle p,q\right\rangle$ satisfying $q\geqslant(p+1)/2$. If $A\!\left(t_r,h/k\right)\gg V$ for $r\in\left\{1,2,\ldots,R\right\}$,
then
\[R\ll k^2\,M^{1+\varepsilon}\,V^{-3}+k^{2q/p}\,M^{1+q/p+\varepsilon}\,V^{-2-(1+2q)/p}.\]
\end{theorem}

The large value result can be turned into a general moment estimate:
\begin{theorem}\label{moments}
Let $M\in\left[1,\infty\right[$,
fix an exponent pair $\left\langle p,q\right\rangle$ satisfying $q\geqslant(p+1)/2$, and let $h$ and $k$ be coprime integers with
$1\leqslant k\ll M^{1/2-\varepsilon}$.
Furthermore, let $\alpha,\beta,\gamma,\delta,A\in\left[0,\infty\right[$ be exponents so that
\[\sum_{n\leqslant x}a(n)\,e\!\left(\frac{nh}k\right)\ll k^\alpha\,x^{\beta+\varepsilon}\]
for sufficiently large $x$ and for $k$ satisfying $x^\gamma\ll k\ll x^\delta$.
Then, for $M^\gamma\ll k\ll M^\delta$,
\[\int\limits_M^{2M}\left|A\!\left(x,\frac hk\right)\right|^A\,\mathrm dx
\ll k^{A/2}\,M^{A/4+1}+\Phi+\Psi,\]
where
\[\Phi=\left\{\begin{array}{ll}
k^{\alpha A+2(1-\alpha)}\,M^{\beta A+(1-2\beta)+\varepsilon}&\text{if $A\geqslant2$},\\
k^{A/2+1}\,M^{A/4+1/2+\varepsilon}&\text{if $A\leqslant2$},
\end{array}\right.\]
and
\[\Psi=\left\{\begin{array}{ll}
k^{\alpha A-\alpha-\alpha/p+(1-\alpha)2q/p}\,M^{\beta A+1-\beta-\beta/p+(1-2\beta)q/p+\varepsilon}&
\text{if $A\geqslant1+\frac{1+2q}p$,}\\[1mm]
k^{A/2-1/2-1/(2p)+q/p}\,M^{A/4+3/4-1/(4p)+q/(2p)+\varepsilon}&\text{if $A\leqslant1+\frac{1+2q}p$.}
\end{array}\right.\]
\end{theorem}

This theorem of course admits many special cases. We will mention only a few.
Using the classical upper bound $A(x,h/k)\ll k^{2/3}\,x^{1/3+\varepsilon}$, certainly valid for $1\ll k\ll M^{1/2-\varepsilon}$, the expressions for $\Phi$ and $\Psi$ simplify to
\[\Phi=\left\{\begin{array}{ll}
k^{2A/3+2/3}\,M^{A/3+1/3+\varepsilon}&\text{if $A\geqslant2$},\\
k^{A/2+1}\,M^{A/4+1/2+\varepsilon}&\text{if $A\leqslant2$},
\end{array}\right.\]
and
\[\Psi=\left\{\begin{array}{ll}
k^{2A/3-2/3-2/(3p)+2q/(3p)}\,M^{A/3+2/3+q/(3p)-1/(3p)+\varepsilon}&
\text{if $A\geqslant1+\frac{1+2q}p$,}\\[1mm]
k^{A/2-1/2-1/(2p)+q/p}\,M^{A/4+3/4-1/(4p)+q/(2p)+\varepsilon}&\text{if $A\leqslant1+\frac{1+2q}p$.}
\end{array}\right.\]

When $M^{1/10}\ll k\ll M^{1/4}$, we can use the better exponents $\alpha=1/4$ and $\beta=3/8$ coming from Theorem~\ref{pointwise-estimate}, and $\Phi$ and $\Psi$ are given by
\[\Phi=\left\{\begin{array}{ll}
k^{A/4+3/2}\,M^{3A/8+1/4+\varepsilon}&\text{if $A\geqslant2$},\\
k^{A/2+1}\,M^{A/4+1/2+\varepsilon}&\text{if $A\leqslant2$},
\end{array}\right.\]
and
\[\Psi=\left\{\begin{array}{ll}
k^{A/4-1/4-1/(4p)+3q/(2p)}\,M^{3A/8+5/8-3/(8p)+q/(4p)+\varepsilon}&
\text{if $A\geqslant1+\frac{1+2q}p$,}\\[1mm]
k^{A/2-1/2-1/(2p)+q/p}\,M^{A/4+3/4-1/(4p)+q/(2p)+\varepsilon}&\text{if $A\leqslant1+\frac{1+2q}p$.}
\end{array}\right.\]

The third moment case $A=3$ is of particular interest to us since we shall use it in the proof of Theorem \ref{fourth-moment-estimate}. Following \cite{Ivic1983}, we choose the exponent pair $\left\langle p,q\right\rangle=\left\langle4/18,11/18\right\rangle=BABAAB\left\langle0,1\right\rangle$, which is sufficient for our purposes. Many other exponent pairs can be found in Table 1 on page 58 in \cite{Montgomery1994}. For this pair $1+\frac{1+2q}p=11$.

In the range $k\ll M^{5/14}$ we then obtain
\[\Phi=k^{2+2/3}\,M^{4/3+\varepsilon}
=k^{8/3}\,M^{4/3+\varepsilon}
\ll k^{3/2}\,M^{7/4+\varepsilon},\]
and
\[\Psi=k^{3/2}\,M^{7/4+\varepsilon}.\]
Thus,
we obtain
\begin{corollary}\label{cubic-moment}
For $k\ll M^{5/14}$, we have
\[\int\limits_M^{2M}\left|A\!\left(x,\frac hk\right)\right|^3\mathrm dx
\ll k^{3/2}\,M^{7/4+\varepsilon}.\]
\end{corollary}
\noindent
By H\"older's inequality and Theorem \ref{mean-square-of-long-sums} below this cannot be improved except possibly for the range of $k$ and the $\varepsilon$.

Because of the $\varepsilon$ in the exponents, these moment estimates are not quite strong enough that we could see that the oscillation result Theorem \ref{oscillation-result} holds for a positive proportion of time. Thus, we need a higher moment estimate where the $\varepsilon$ has been shaved away. For this purpose, as well as for its intrinsic interest, we derive a reasonably sharp fourth moment estimate.
\begin{theorem}\label{fourth-moment-estimate}
For $M\in\left[1,\infty\right[$ and for coprime integers $h$ and $k$ with $1\leqslant k\ll M^{1/6-\varepsilon}$, we have
\[\int\limits_M^{2M}\left|A\!\left(x,\frac hk\right)\right|^4\mathrm dx
=C_F\,k^2\,M^2+O(k^{11/4}\,M^{15/8+\varepsilon})+O(k^{13/6}\,M^{23/12+\varepsilon}),\]
where the constant coefficient $C_F$ is
\begin{multline*}
\frac3{64\pi^4}\left(\,\sum_{\substack{a,b,c,d=1\\\sqrt{\vphantom ba}+\sqrt b=\sqrt{\vphantom dc}+\sqrt d}}^\infty
+\sum_{\substack{a,b,c,d=1\\\sqrt{\vphantom ba}+\sqrt{\vphantom dc}=\sqrt b+\sqrt d}}^\infty
+\sum_{\substack{a,b,c,d=1\\\sqrt{\vphantom ba}+\sqrt d=\sqrt b+\sqrt{\vphantom dc}}}^\infty\,\right)
\left(abcd\right)^{-3/4}\\
\cdot a(a)\,a(b)\,\overline{a(c)\,a(d)}\,e\!\left(\frac{-(a+b-c-d)\,\overline h}k\right).
\end{multline*}
\end{theorem}
\noindent
The constant $C_F$ is a finite positive real number by Deligne's estimate $a(n)\ll d(n)$ and the estimates (3.6) in \cite{Tsang1992}. We remark that Theorem \ref{fourth-moment-estimate} immediately implies that Corollary \ref{cubic-moment} holds without the $\varepsilon$ in the upper bound when $k\ll M^{1/6-\varepsilon}$.

Next, we turn to the mean square estimates for short linear sums.
\begin{theorem}\label{mean-square-of-very-short-sums}
Let $M\gg1$, let $1\leqslant\Delta\ll\sqrt M$, let $h$ and $k$ be coprime integers with $1\leqslant k\ll\Delta^{1/2}\,M^{-\varepsilon}$, and, finally, let $k^2\,M^{1+\varepsilon}\,\Delta^{-1}\ll\Xi\ll M$. Then
\[\int\limits_M^{M+\Xi}\,\left|\sum_{x\leqslant n\leqslant x+\Delta}a(n)\,e\!\left(\frac{nh}k\right)\right|^2\mathrm dx\ll\Xi\,\Delta.\]
Furthermore, if $\Delta\ll M^{1/2-\varepsilon}$, and $M\gg_\varepsilon1$, and if the underlying cusp form does not vanish identically, then the second moment is $\asymp\Xi\,\Delta$.
\end{theorem}
\noindent 
Next, a slightly different upper bound which treats sums of many different lengths at once.
\begin{theorem}\label{sharper-mean-square-for-short-sums}
Let $M\gg1$, let $M^\varepsilon\ll\Delta\ll\sqrt M$, let $h$ and $k$ be coprime integers with $1\leqslant k\ll\Delta^{1/4}\,M^{-\varepsilon}$.
Then
\[\int\limits_M^{2M}\max_{0\leqslant U\leqslant\Delta}
\left|\sum_{x\leqslant n\leqslant x+U}a(n)\,e\!\left(\frac{nh}k\right)\right|^2\mathrm dx\ll M\,\Delta\,\log^2M.\]
\end{theorem}

Finally, we arrive at our final theorem:
\begin{theorem}\label{oscillation-result}
Let $M\in\mathbb R_+$ be sufficiently large, and let $h$ and $k$ be coprime integers with $1\leqslant k\ll M^{1/8-\varepsilon}$. Then there are $\gg M^{1/2}\,\log^2M$ pairwise disjoint subintervals $I$ of $\left[M,2M\right]$ of length $\gg M^{1/2}\log^{-2}M$ such that, for $x\in I$,
\[\sum_{n\leqslant x}a(n)\,e\!\left(\frac{nh}k\right)
\gg k^{1/2}\,M^{1/4},\]
and we may even demand that the argument of the above sum as a complex number belongs to a given interval of length $3\pi/2+\varepsilon$ for $x\in I$.
\end{theorem}
\noindent
Here, of course, the intervals will depend on $M$, $h$ and $k$, but the implicit constants will not.

\section{Various lemmas and well-known theorems}

We collect in this section the pre-existing results and some of their consequences required in the proofs.

\subsection{Results on holomorphic cusp forms}\label{results-on-holomorphic-cusp-forms}

Deligne's famous work on Weil's conjectures gives as an application \cite{Deligne1974} an estimate for individual Fourier coefficients of a holomorphic cusp form:
\[a(n)\ll d(n),\]
and in fact, $\left|a(n)\right|\leqslant d(n)$, if the cusp form in question is an arithmetically normalized Hecke eigenform. This is often combined with Shiu's estimate \cite{Shiu1980} for the divisor function:
\[\sum_{x\leqslant n\leqslant x+y}d(n)\ll y\,\log x,\]
for $x^\varepsilon\ll y\ll x$.

Rankin \cite{Rankin1939} and Selberg \cite{Selberg1940} obtained the main term for the mean square of the Fourier coefficients:
\begin{theorem}
We have
\[\sum_{n\leqslant x}\left|a(n)\right|^2=Cx+O\bigl(x^{3/5}\bigr),\]
where $C\in\mathbb R_+$ only depends on the cusp form in question.
\end{theorem}
Through integration by parts, this can be used in many forms. For example, for $x\in\left[1,\infty\right[$, $\alpha\in\left]0,1\right[$ and $\beta\in\left]1,\infty\right[$, we have
\[\sum_{n\leqslant x}\frac{\left|a(n)\right|^2}{n^\alpha}
\ll_\alpha x^{1-\alpha}\quad\text{and}\quad
\sum_{x>n}\frac{\left|a(n)\right|^2}{n^\beta}
\ll_\beta x^{1-\beta}.\]

Most of our proofs depend on the important truncated Voronoi identity, a good presentation of which is given in Chapter 1 of Jutila's monograph \cite{Jutila1987a}, where this identity appears in Theorem 1.1.
\begin{theorem}
Let $x\in\left[1,\infty\right[$, let $N\in\left[1,\infty\right[$, assume that $N\ll x$, let $h$ be an integer, let $k$ be a positive integer, and assume that $\left(h,k\right)=1$ and $k\leqslant x$. Then
\begin{multline*}
\sum_{n\leqslant x}a(n)\,e\!\left(\frac{nh}k\right)
=\frac{k^{1/2}\,x^{1/4}}{\pi\,\sqrt2}
\sum_{n\leqslant N}a(n)\,n^{-3/4}\,e\!\left(\frac{-n\overline h}k\right)\\\cdot\cos\!\left(\frac{4\pi\sqrt{nx}}k-\frac\pi4\right)
+O\bigl(k\,x^{1/2+\varepsilon}\,N^{-1/2}\bigr).
\end{multline*}
\end{theorem}
\noindent In Theorem 1.2 of the same chapter, Jutila also gives as an application the mean square estimate for long linear sums, which with our normalization for $a(n)$ reads as follows.
\begin{theorem}\label{mean-square-of-long-sums}
Let $M\geqslant1$, and let $h$ and $k$ be coprime integers with $1\leqslant k\ll M$. Then
\[\int\limits_1^M\left|A\!\left(x,\frac hk\right)\right|^2\mathrm dx
=C\,k\,M^{3/2}
+O\bigr(k^2\,M^{1+\varepsilon}\bigr)
+O\bigl(k^{3/2}\,M^{5/4+\varepsilon}\bigr),\]
where
\[C=\frac1{6\pi^2}\sum_{n=1}^\infty\frac{\left|a(n)\right|^2}{n^{3/2}}.\]
\end{theorem}
\noindent The proof of Theorem \ref{pointwise-estimate} requires the full Voronoi type summation formula for holomorphic cusp forms.
\begin{theorem}\label{voronoi-summation-formula}
For positive real numbers $a$ and $b$ with $a<b$, and for coprime integers $h$ and $k$ with $k\geqslant1$, and a continuously differentiable function $f\colon\mathbb R_+\longrightarrow\mathbb R_+$, we have
\begin{multline*}
\sum_{a\leqslant n\leqslant b}'a(n)\,e\!\left(\frac{nh}k\right)f(n)\\
=(-1)^{\kappa/2}\,\frac{2\pi}k\sum_{n=1}^\infty
a(n)\,e\!\left(\frac{-n\overline h}k\right)
\int\limits_a^bf(x)\,J_{\kappa-1}\!\left(\frac{4\pi\sqrt{nx}}k\right)\mathrm dx.
\end{multline*}
Here $\sum'$ means that if $a$ (or $b$) is an integer, then the term $f(a)$ (or $f(b)$) should be halved.
\end{theorem}
\noindent
Again, the presentation of Chapter 1 of \cite{Jutila1987a} is recommended; Theorem 1.7 there gives the above summation formula.

\subsection{Tools for exponential sums without Fourier coefficients}

The following lemma due to Bombieri, which is Lemma 1.5 in \cite{Montgomery1971}, is not strictly speaking an exponential sum result, but we will use it in the proof of Theorem~\ref{large-values} to separate the Fourier coefficients from the exponential sums, leaving only plain exponential sums which then can be estimated using the theory of exponent pairs.
\begin{theorem}\label{bombieri-lemma}
Let $H$ be a Hilbert space. Denote its inner product by $\left\langle\cdot\middle|\cdot\right\rangle$ and its norm by $\left\|\cdot\right\|$. Also, let $R\in\mathbb Z_+$, and let $\xi$, $\varphi_1$, $\varphi_2$, \dots, $\varphi_R$ be vectors of $H$. Then
\[\sum_{r=1}^R\left|\left\langle\xi\middle|\varphi_r\right\rangle\right|^2
\leqslant\left\|\xi\right\|^2\max_{1\leqslant r\leqslant R}\sum_{s=1}^R\left|\left\langle\varphi_r\middle|\varphi_s\right\rangle\right|.\]
\end{theorem}
We will actually only apply this for $H=\mathbb C^N$ for some $N\in\mathbb Z_+$, and the inner product will be the usual one so that for $z=\left\langle z_1,\ldots,z_N\right\rangle\in\mathbb C^N$ and $w=\left\langle w_1,\ldots,w_N\right\rangle\in\mathbb C^N$ we have
\[\left\langle z\middle|w\right\rangle=\overline z_1w_1+\ldots+\overline z_Nw_N,\quad\text{and}\quad
\left\|z\right\|^2=\left|z_1\right|^2+\ldots+\left|z_N\right|^2.\]

An excellent presentation of the theory of exponent pairs can be found in~\cite{GrahamEtKolesnik1991}. Other presentations on the topic can be found in e.g. Chapter 2 of \cite{Ivic2003} and in Chapter 3 of \cite{Montgomery1994}. For the purposes of the proof of Theorem \ref{large-values}, it is enough to state here that if the pair $\left\langle p,q\right\rangle\in\left[0,\frac12\right]\times\left[\frac12,1\right]$ happens to be an exponent pair, then we may estimate
\[\sum_{M\leqslant n\leqslant M+\Delta}e(A\sqrt n)\ll
A^p\,M^{q-p/2}+A^{-1}\,M^{1/2},\]
for positive real constants $A$. The square root $\sqrt n$ will appear from the right-hand side of the truncated Voronoi identity.

\subsection{Lemmas on exponential integrals}

\begin{lemma}\label{primitive-lemma}
Let $M\geqslant1$, let $0\leqslant\Xi\leqslant M$, and let $k$ and $n$ be positive integers. Then
\[\int\limits_M^{M+\Xi}x^{1/4}\,\cos\!\left(\frac{4\pi\sqrt{nx}}k-\frac\pi4\right)\mathrm dx\ll\frac k{\sqrt n}\,M^{3/4}.\]
\end{lemma}
\noindent This is a typical consequence of the first derivative test, which is discussed e.g. in Section 2.1 of \cite{Ivic2003}. More precisely, the claim follows from Lemma 2.1 of \cite{Ivic2003} and integration by parts.

\begin{lemma}\label{mean-square-lemma-1}
Let $M\gg1$, $0\leqslant U\ll\sqrt M$, $0\leqslant T\ll\sqrt M$, $0\leqslant\Xi\ll M$, and let $m$, $n$ and $k$ be positive integers. Then
\[\int\limits_M^{M+\Xi}x^{1/2}\,e\!\left(\pm\!\left(\frac{2\sqrt{m\left(x+T\right)}}k+\frac{2\sqrt{n\left(x+U\right)}}k\right)\!\right)\mathrm dx
\ll\frac{k\,M}{\sqrt m+\sqrt n}.\]
\end{lemma}
\noindent This follows immediately from the first derivative test.

\begin{lemma}\label{mean-square-lemma-2}
Let $M\gg1$, $0\leqslant U\leqslant T\ll\sqrt M$, $0\leqslant\Xi\ll M$, and let $m$, $n$ and $k$ be positive integers with $m<n$. Then
\[\int\limits_M^{M+\Xi}x^{1/2}\,e\!\left(\pm\!\left(\frac{2\sqrt{m\left(x+T\right)}}k-\frac{2\sqrt{n\left(x+U\right)}}k\right)\!\right)\mathrm dx
\ll\frac{k\,\sqrt{n}\,M}{n-m}.\]
\end{lemma}
\noindent This also follows immediately from the first derivative test given the simple observation that
\[\frac{k\,\sqrt M}{\sqrt n-\sqrt m}
\ll\frac{k\,\sqrt{n\,M}}{n-m}.\]

\begin{lemma}\label{mean-square-lemma-3}
Let $M\gg1$, $0<T\ll\sqrt M$, $0\leqslant\Xi\ll M$, and let $m$, $n$ and $k$ be positive integers with $m<n$ and $m\leqslant\frac M{2T}$. Then
\[\int\limits_M^{M+\Xi}x^{1/2}\,e\!\left(\pm\!\left(\frac{2\sqrt{mx}}k-\frac{2\sqrt{n\left(x+T\right)}}k\right)\!\right)\mathrm dx
\ll\frac{k\,\sqrt{n}\,M}{n-m}.\]
\end{lemma}
\noindent This is once more a consequence of the first derivative test. The relevant estimate is
\begin{align*}
&\frac{\mathrm d}{\mathrm dx}\!\left(\frac{2\sqrt{n\left(x+T\right)}}k-\frac{2\sqrt{mx}}k\right)
=\frac{\sqrt n}{k\,\sqrt{x+T}}-\frac{\sqrt m}{k\,\sqrt x}\\
&=\frac{nx-mx-mT}{k\,\sqrt{x\left(x+T\right)}\left(\sqrt{nx}+\sqrt{m\left(x+T\right)}\right)}
\gg\frac{\left(n-m\right)M}{k\,M\left(\sqrt n+\sqrt m\right)\sqrt M}
\gg\frac{n-m}{k\,\sqrt{nM}}.
\end{align*}

Finally, one more corollary of the first derivative test:
\begin{lemma}\label{fourth-moment-lemma}
Let $M\in\left[1,\infty\right[$, let $k\in\mathbb Z_+$, and let $\Delta\in\mathbb R_+$. Then
\[\int\limits_M^{2M}x\,e\!\left(\pm\frac{\Delta\,\sqrt{x}}k\right)\mathrm dx\ll\frac{k\,M^{3/2}}\Delta.\]
\end{lemma}

\subsection{Lemmas on the spacing of square roots}

The following three lemmas are Lemmas 2, 5 and 6 in \cite{IvicEtSargos2007}, and they will be used in the consideration of the fourth moment estimate in Section \ref{fourth-moment-section}.

\begin{lemma}\label{lower-bound-for-square-root-expression}
If $a$, $b$, $c$ and $d$ are positive integers with $a\leqslant c$, $b\leqslant c$ and $\sqrt{\vphantom da}+\sqrt b\pm\sqrt{\vphantom bc}\neq\sqrt d$, then
\[\bigl|\sqrt{\vphantom da}+\sqrt b\pm\sqrt{\vphantom bc}-\sqrt d\bigr|
\gg c^{-2}\,(abc)^{-1/2}.\]
\end{lemma}

\begin{lemma}\label{counting-two-two-quadruples}
Let $A$, $B$, $C$, $D$ and $\delta$ be positive real numbers with $1\ll A\ll C$, $1\ll B\ll C$ and $1\ll D\ll C$. Then the number of quadruples of positive integers $a$, $b$, $c$ and $d$ with $A<a\leqslant2A$, $B<b\leqslant2B$, $C<c\leqslant2C$, $D<d\leqslant2D$ and
\[\bigl|\sqrt{\vphantom da}+\sqrt b-\sqrt{\vphantom dc}-\sqrt d\,\bigr|\leqslant\delta\sqrt C\]
is both
\[\ll ABCD\bigl(\delta+C^{-3/2+\varepsilon}\bigr)
+\bigl(ABCD\bigr)^{1/2}\]
and
\[\ll ABCD\bigl(\delta C^2+\left(ABCD\right)^{-1/2}\bigr)C^\varepsilon.\]
\end{lemma}

\begin{lemma}\label{counting-three-one-quadruples}
Let $A$, $B$, $C$, $D$ and $\delta$ be positive real numbers with $1\ll A\ll C$ and $1\ll B\ll C$. Then the number of quadruples of positive integers $a$, $b$, $c$ and $d$ with $A<a\leqslant2A$, $B<b\leqslant2B$, $C<c\leqslant2C$, $D<d\leqslant2D$ and
\[0<\bigl|\sqrt{\vphantom da}+\sqrt b+\sqrt{\vphantom dc}-\sqrt d\,\bigr|\leqslant\delta\sqrt C\]
is both
\[\ll ABCD\bigl(\delta+C^{-3/2+\varepsilon}\bigr)\]
and
\[\ll ABCD\bigl(\delta C^2+\left(ABCD\right)^{-1/2}\bigr)C^\varepsilon.\]
\end{lemma}

\section{Proof of Theorem \ref{pointwise-estimate}}

Theorem \ref{pointwise-estimate} will follow by introducing a smooth weight function and using the following weighted result:
\begin{lemma}\label{weighted-pointwise-estimate}
Let $M\in\left[1,\infty\right[$, let $\Delta\in\left[1,M\right]$, let $h$ and $k$ be coprime integers with $1\leqslant k\ll M^{1/2}$, and let $X\in\left[1,\infty\right[$. Also, let $U\in\left[1,\Delta/3\right]$ and let $w\in C^2(\mathbb R)$ be supported in $\left[M,M+\Delta\right]$ and satisfy
$w\ll1$, $w'\ll U^{-1}$ and $w''\ll U^{-2}$, as well as $w(x)=1$ for $x\in\left[M+U,M+\Delta-U\right]$. Then
\[\sum_{M\leqslant n\leqslant M+\Delta}a(n)\,e\!\left(\frac{nh}k\right)w(n)
\ll k^{1/2}\,M^{1/4}\,X^{1/4}+k^{3/2}\,X^{-1/4}\,U^{-1}\,M^{3/4}.\]
\end{lemma}

Remark: By examining the proof of the lemma, one sees that the first term can be omitted by choosing $X=1/2$ with the upper bound then becoming $\ll k^{3/2}\,U^{-1}\,M^{3/4}$. For example, if $U\gg k^{3/2}\,M^{3/4+\varepsilon}$, then the weighted sum in question is $o(1)$ as $M\longrightarrow\infty$.

\paragraph{Proof of Theorem \ref{pointwise-estimate}.} The last estimate of Theorem \ref{pointwise-estimate} follows from the results in \cite{Jutila1987b}. The other estimates follow by first smoothing the sum with a weight function $w$ satisfying the assumptions of Lemma \ref{weighted-pointwise-estimate}, and then applying Lemma \ref{weighted-pointwise-estimate} to the smooth sum. However, the different ranges of $k$ will be handled slightly differently.

First, if $M^{1/10}\ll k\ll M^{1/4}$, then we shall choose $U=k^{3/2}\,M^{1/4}$. It is easy to check that $M^{1/10}\ll k\ll M^{1/4}$ implies that $M^{2/5}\ll U\ll M^{5/8}$, so that Theorem 5.5 from \cite{ErnvallHytonenEtKarppinen2008} applies and the smoothing error will be $\ll U^{1/6}\,M^{1/3+\varepsilon}$ by a simple integration by parts argument. Combining this with Lemma \ref{weighted-pointwise-estimate} gives
\[\sum_{M\leqslant n\leqslant 2M}a(n)\,e\!\left(\frac{nh}k\right)
\ll U^{1/6}\,M^{1/3+\varepsilon}+k^{1/2}\,M^{1/4}\,X^{1/4}+k^{3/2}\,X^{-1/4}\,U^{-1}\,M^{3/4}.\]
Choosing now $X=k^2\,M\,U^{-2}$ this simplifies to $\ll k^{1/4}\,M^{3/8+\varepsilon}$, as required.

Similarly, if $M^{1/4}\ll k\ll M^{19/64}$, we choose $U=k^{2/3}\,M^{11/24}$ which is $\gg M^{5/8}$ and $\ll M^{11/16}$, so that Theorem 5.7 from \cite{ErnvallHytonenEtKarppinen2008} gives the smoothing error $\ll U\,M^{-3/16+\varepsilon}$. Choosing then $X=k^{2/3}\,M^{1/12}$ in Lemma \ref{weighted-pointwise-estimate} gives the upper bound $\ll k^{2/3}\,M^{13/48+\varepsilon}$.

In the same vein, if $M^{19/64}\ll k\ll M^{21/64}$, then we may choose $U=M^{23/32}$ and use Theorem 5.16 from \cite{ErnvallHytonenEtKarppinen2008} to get the smoothing error $\ll M^{15/32+\varepsilon}$. Then choosing $X=k^{-2}\,M^{7/8}$ in Lemma \ref{weighted-pointwise-estimate} gives only terms which are $\ll M^{15/32+\varepsilon}$.

In the range $M^{21/64}\ll k\ll M^{3/8}$ we choose $U=k^{2/3}\,M^{1/2}$ and apply Theorem 5.16 from \cite{ErnvallHytonenEtKarppinen2008} to get the smoothing error $\ll U\,M^{-1/4+\varepsilon}$. Then, choosing $X=k^{2/3}$ in Lemma \ref{weighted-pointwise-estimate} gives the upper bound $\ll k^{2/3}\,M^{1/4+\varepsilon}$.

Finally, for $k\ll\sqrt M$, we may choose $U=k^{2/3}\,M^{1/3}\log^{-2/3}M$ and $X=k^{2/3}\,M^{1/3}\,\log^{4/3}M$, and the smoothing error can be estimated by Deligne's estimate and Shiu's theorem by $\ll k^{2/3}\,M^{1/3}\,\log^{1/3}M$, which is also the upper bound that Lemma \ref{weighted-pointwise-estimate} gives for the smooth sum.

\paragraph{Proof of Lemma \ref{weighted-pointwise-estimate}.}
An application of the additively twisted Voronoi-type summation formula for holomorphic cusp forms gives
\begin{multline*}\sum_{M\leqslant n\leqslant M+\Delta}a(n)\,e\!\left(\frac{nh}k\right)w(n)\\
=(-1)^{\kappa/2}\frac{2\pi}k\sum_{n=1}^\infty
a(n)\,e\!\left(\frac{-n\overline h}k\right)
\int\limits_M^{M+\Delta}J_{\kappa-1}\!\left(\frac{4\pi\sqrt{nx}}k\right)
w(x)\,\mathrm dx.\end{multline*}
We shall consider the series in two parts $\sum_{n\leqslant X}+\sum_{n>X}$.

In the terms with $n\leqslant X$ we first apply the asymptotics of the Bessel $J$-function
\[J_{\nu}\!\left(x\right)=\sqrt{\frac2{\pi x}}\,
\cos\!\left(x-\frac{\nu\pi}2-\frac\pi4\right)
+O_\nu\!\left(x^{-3/2}\right),\]
and then apply the first derivative test on the integrals involving the main terms of the $J$-asymptotics and estimates by absolute values on the integrals involving the error term of the $J$-asymptotics:
\begin{align*}
&\int\limits_M^{M+\Delta}J_{\kappa-1}\!\left(\frac{4\pi\sqrt{nx}}k\right)
w\!\left(x\right)\mathrm dx\\
&\qquad=\frac{k^{1/2}}{\pi\,\sqrt2}\,n^{-1/4}\int\limits_M^{M+\Delta}
x^{-1/4}\,\cos\!\left(\frac{4\pi\sqrt{nx}}k-\frac{\kappa\pi}2+\frac\pi4\right)w(x)\,\mathrm dx\\
&\qquad\qquad+O\!\left(\Delta\,k^{3/2}\,n^{-3/4}\,M^{-3/4}\right).\\
&\qquad\ll k^{1/2}\,n^{-1/4}\,M^{-1/4}\,k\,n^{-1/2}\,M^{1/2}
+\Delta\,k^{3/2}\,n^{-3/4}\,M^{-3/4}\\
&\qquad\ll k^{3/2}\,n^{-3/4}\,M^{1/4}.
\end{align*}
Now the sum over $n\leqslant X$, using the Rankin--Selberg estimate with the Cauchy--Schwarz inequality, gives the first term
$\ll k^{1/2}\,X^{1/4}\,M^{1/4}$.

Next, we treat the terms with $n>X$ by integrating by parts twice using the fact that
\[\frac{\mathrm d}{\mathrm dx}\left(x^\nu\,J_\nu(x)\right)=x^\nu\,J_{\nu-1}(x),\]
and estimating by absolute values,
and get
\begin{align*}
&\int\limits_M^{M+\Delta}J_{\kappa-1}\!\left(\frac{4\pi\sqrt{nx}}k\right)w(x)\,\mathrm dx\\
&\quad=\frac{k^2}{4\pi^2 n}\int\limits_M^{M+\Delta}
J_{\kappa+1}\!\left(\frac{4\pi\sqrt{nx}}k\right)
\left(x\,w''(x)+\left(1-\kappa\right)w'(x)+\frac{\kappa^2-1}{4x}\cdot w(x)\right)\mathrm dx\\
&\quad\ll k^{5/2}\,n^{-5/4}\,M^{3/4}\,U^{-1}.
\end{align*}
It is important to observe here that we really get just $U^{-1}$ instead of $\Delta\,U^{-2}$, the reason being that after integrations by parts the terms involving $w'$ and $w''$ in the resulting integrand vanish outside the support of $w'(x)$ which is of length $\ll U$.
Thus the series over $n>X$ gives the second term $\ll k^{3/2}\,X^{-1/4}\,M^{3/4}\,U^{-1}$.

\section{Proof of Theorem \ref{large-values}}

Let us restrict to values of $M$ so large that $k\leqslant M$, and let $N\in\left[1,\infty\right[$ be a parameter such that $1\leqslant N\ll M$. We shall choose the precise value of the parameter $N$ later.
From the truncated Voronoi identity, i.e.\ Theorem \ref{truncated-voronoi-identity} above, we now see that for $x\in\left[M,2M\right]$,
\begin{align*}
&\sum_{n\leqslant x}a(n)\,e\!\left(\frac{nh}k\right)\\
&\qquad=\frac{k^{1/2}\,x^{1/4}}{\pi\,\sqrt2}\sum_{n\leqslant N}a(n)\,n^{-3/4}\,e\!\left(\frac{-n\overline h}k\right)\,\cos\left(\frac{4\pi\sqrt{nx}}k-\frac\pi4\right)\\
&\qquad\qquad+O\!\left(k\,x^{1/2+\varepsilon}\,N^{-1/2}\right).\\
&\qquad\ll\sum_\pm k^{1/2}\,x^{1/4}\left|\sum_{n\leqslant N}a(n)\,n^{-3/4}\,e\!\left(\frac{-n\overline h}k\right)\,e\!\left(\pm\frac{2\sqrt{nx}}k\right)\right|\\
&\qquad\qquad+k\,x^{1/2+\varepsilon}\,N^{-1/2}.
\end{align*}
Here the sign $\sum_\pm$ means summing over both choices of the sign $\pm$.

We shall divide the interval $\left[M,2M\right]$ into subintervals of length not exceeding $M_0\in\mathbb R_+$. The value of $M_0$ will be chosen near the end of the argument. We shall consider the number $R_0$ of the points $t_1$, \dots, $t_R$ lying in a single subinterval, and we shall call these points $t_1$, $t_2$, \dots, $t_{R_0}$. Once an upper bound $R_0\ll\Upsilon$ has been obtained, in which $\Upsilon$ does not depend on the subinterval, the estimate for $R$ will be obtained from that of $R_0$ via
\[R\ll\Upsilon\left(1+\frac M{M_0}\right).\]

Now, applying the above inequality $R_0$ times gives
\begin{align*}
R_0\,V^2
&\ll\sum_\pm\sum_{r\leqslant R_0}\left|k^{1/2}\,M^{1/4}\sum_{n\leqslant N}a(n)\,n^{-3/4}\,e\!\left(\frac{-n\overline h}k\right)\,e\!\left(\pm\frac{2\sqrt{nt_r}}k\right)\right|^2\\
&\qquad+k^2\,R_0\,M^{1+\varepsilon}\,N^{-1}.
\end{align*}
The error term coming from the truncated Voronoi identity may be absorbed to the left-hand side provided that $N\gg k^2\,V^{-2}\,M^{1+\varepsilon}$, and we will actually chooce $N=k^2\,M^{1+\varepsilon}\,V^{-2}$, where the $\varepsilon$ is so small that the condition $k\,M^\varepsilon\ll V$ implies $N\ll M$.

Now we are ready to continue our estimates by
\begin{align*}
&R_0\,V^2\ll\sum_\pm\sum_{r\leqslant R_0}\left|k^{1/2}\,M^{1/4}\sum_{n\leqslant N}a(n)\,n^{-3/4}\,e\!\left(\frac{-n\overline h}k\right)\,e\!\left(\pm\frac{2\sqrt{nt_r}}k\right)\right|^2\\
&\ll k\,M^{1/2}\sum_\pm\sum_{r\leqslant R_0}\left|\sum_{\substack{U\leqslant N/2\\\text{dyadic}}}\sum_{U<n\leqslant2U}a(n)\,n^{-3/4}\,e\!\left(\frac{-n\overline h}k\right)\,e\!\left(\pm\frac{2\sqrt{nt_r}}k\right)\right|^2\\
&\ll k\,M^{1/2}\,\log N\sum_\pm\sum_{r\leqslant R_0}\sum_{\substack{U\leqslant N/2\\\text{dyadic}}}\left|\sum_{U<n\leqslant2U}a(n)\,n^{-3/4}\,e\!\left(\frac{-n\overline h}k\right)\,e\!\left(\pm\frac{2\sqrt{nt_r}}k\right)\right|^2\\
&\ll k\,M^{1/2}\,\log^2N\sum_\pm\max_{U\leqslant N/2}\sum_{r\leqslant R_0}\left|\sum_{U<n\leqslant2U}a(n)\,n^{-3/4}\,e\!\left(\frac{-n\overline h}k\right)\,e\!\left(\pm\frac{2\sqrt{nt_r}}k\right)\right|^2\!.
\end{align*}
We now apply Bombieri's lemma, i.e.\ Theorem \ref{bombieri-lemma} above, to get
\begin{align*}
R_0\,V^2&\ll k\,M^{1/2+\varepsilon}\,\max_{U\leqslant N/2}U^{-1/2}\\
&\qquad\cdot\max_{r\leqslant R_0}\sum_{s\leqslant R_0}\left|\sum_{U<n\leqslant2U}e\!\left(\frac{2\sqrt n\left(\sqrt{t_r}-\sqrt{t_s}\right)}k\right)\right|\\
&\ll k\,M^{1/2+\varepsilon}\,N^{1/2}\\
&\qquad+k\,M^{1/2+\varepsilon}\,\max_{U\leqslant N/2}U^{-1/2}\\
&\qquad\qquad\cdot\max_{r\leqslant R_0}\sum_{\substack{s\leqslant R_0\\s\neq r}}\left|\sum_{U<n\leqslant2U}e\!\left(\frac{2\sqrt n\left(\sqrt{t_r}-\sqrt{t_s}\right)}k\right)\right|.
\end{align*}

Next we get to apply the exponent pair $\left\langle p,q\right\rangle$. We have
\[\left|\sqrt{t_r}-\sqrt{t_s}\right|\asymp\left|\,\int\limits_{t_r}^{t_s}\frac{\mathrm dt}{t^{1/2}}\right|\asymp\frac{\left|t_r-t_s\right|}{M^{1/2}}.\]
This is both $\ll M_0\,M^{-1/2}$ and $\gg \left|t_r-t_s\right|\,M^{-1/2}$. Therefore we have the estimate
\[\sum_{U<n\leqslant2U}e\!\left(\frac{2\sqrt n\left(\sqrt{t_r}-\sqrt{t_s}\right)}k\right)\ll k^{-p}\,M_0^p\,M^{-p/2}\,U^{q-p/2}+\frac{k\,U^{1/2}\,M^{1/2}}{\left|t_r-t_s\right|}.\]
Substituting this back to the previous estimates gives, estimating $U^{q-p/2-1/2}\ll N^{q-p/2-1/2}$,
\begin{align*}
R_0\,V^2
&\ll k\,M^{1/2+\varepsilon}\,N^{1/2}
+k^2\,M^{1+\varepsilon}\,V^{-1}\\
&\qquad+k^{1-p}\,R_0\,M_0^p\,M^{1/2-p/2+\varepsilon}\,N^{q-p/2-1/2}.
\end{align*}

The choice $N=k^2\,M^{1+\varepsilon}\,V^{-2}$ allows us to merge the first two terms on the right-hand side giving
\[R_0\ll k^2\,M^{1+\varepsilon}\,V^{-3}
+k^{2(q-p)}\,R_0\,M_0^p\,M^{q-p+\varepsilon}\,V^{p-2q-1}.\]
In order to absorb the last term to the left-hand side, the choice of $M_0$ should be such that
\[k^{2(q-p)}\,M_0^p\,M^{q-p+\varepsilon}\,V^{p-2q-1}\ll1,\]
with a sufficiently small implicit constant, of course.
We choose $M_0^p$ to be as large as possible, so that
\[M_0^p\asymp k^{2(p-q)}\,M^{p-q-\varepsilon}\,V^{1+2q-p},\]
and we obtain
\[R_0\ll k^2\,M^{1+\varepsilon}\,V^{-3}.\]
Finally, combining the above considerations leads to the desired estimate
\begin{align*}
R&\ll k^2\,M^{1+\varepsilon}\,V^{-3}
+k^2\,M_0^{-1}\,M^{2+\varepsilon}\,V^{-3}\\
&\ll k^2\,M^{1+\varepsilon}\,V^{-3}
+k^{2q/p}\,M^{1+q/p+\varepsilon}\,V^{-2-(1+2q)/p}.
\end{align*}

\section{Proof of Theorem \ref{moments}}

So, we need to consider the integral
\[\int\limits_{M}^{2M}\left|A\!\left(x,\frac hk\right)\right|^A\,\mathrm dx.\]
The parts where the integrand is $\ll k^{A/2}\,M^{A/4}$ are estimated by absolute values to be $\ll k^{A/2}\,M^{A/4+1}$.
The remaining range of values of $\left|A(x,h/k)\right|$, namely the interval $\left[k^{1/2}\,M^{1/4},k^{\alpha}\,M^{\beta+\varepsilon}\right]$, is divided dyadically into intervals of the shape $\left[V,2V\right]$. This might require extending two of the subintervals, but this doesn't matter. The number of such subintervals is $\ll\log M$. For each interval, we choose the maximum possible number of points, say $R(V)$ points, from the interval $\left[M,2M\right]$, at which the size of $A\!\left(x,h/k\right)$ is in $\left[V,2V\right]$. We choose the points so that they are spaced with distance at least $V$ between each pair of them.

The integral in the regions, where the integrand is large, will be
\begin{align*}
&\ll\sum_VV\cdot R(V)\,V^A\\
&\ll\sum_V\left(k^2\,M^{1+\varepsilon}\,V^{-3}+k^{2q/p}\,M^{1+q/p+\varepsilon}\,V^{-2-(1+2q)/p}\right)V^{A+1}.
\end{align*}

The result now follows easily: the first term in the parantheses gives rise to $\Phi$, the second to $\Psi$. The upper bounds of $\Phi$ and $\Psi$ are then obtained just by estimating $V$ appropriately from below by $k^{1/2}\,M^{1/4}$ or from above by $k^\alpha\,M^{\beta+\varepsilon}$.

\section{Proof of Theorem \ref{fourth-moment-estimate}}\label{fourth-moment-section}

Let $1\leqslant N\ll M$ be arbitrary for the present time. The truncated Voronoi identity, ie.\ Theorem \ref{truncated-voronoi-identity} above, gives
\begin{multline*}
\int\limits_M^{2M}\left|A\!\left(x,\frac hk\right)\right|^4\mathrm dx\\
=\frac{k^2}{4\pi^4}\int\limits_M^{2M}x\left|\sum_{n\leqslant N}a(n)\,n^{-3/4}\,e\!\left(\frac{-n\overline h}k\right)\cos\!\left(\frac{4\pi\sqrt{nx}}k-\frac\pi4\right)\right|^4\mathrm dx
+\text{error},
\end{multline*}
where the error is
\begin{align*}
&\ll\int\limits_M^{2M}\,\left|k^{1/2}\,x^{1/4}\sum_{n\leqslant N}\frac{a(n)}{n^{3/4}}\,e\!\left(\frac{-n\overline h}k\right)\cos\left(\frac{4\pi\sqrt{nx}}k-\frac\pi4\right)\right|^3k\,x^{1/2+\varepsilon}\,N^{-1/2}\,\mathrm dx\\
&\qquad+M\bigl(k\,M^{1/2+\varepsilon}\,N^{-1/2}\bigr)^4\\
&\ll k\,M^{1/2+\varepsilon}\,N^{-1/2}\int\limits_M^{2M}\left|A\!\left(n,\frac hk\right)\right|^3\mathrm dx+M\bigl(k\,M^{1/2+\varepsilon}\,N^{-1/2}\bigr)^4\\
&\ll k^{5/2}\,M^{9/4+\varepsilon}\,N^{-1/2}+k^4\,M^{3+\varepsilon}\,N^{-2}.
\end{align*}

To treat the main term, we will expand the fourth power as $\left|\Sigma\right|^4=\overline\Sigma\vphantom\Sigma^2\,\Sigma^2$, exchange the order of integration and summation, and write the cosines in terms of exponential functions. This, of couse, leads to a large number of terms of the form
\begin{multline*}
\frac{k^2}{64\pi^4}\,a(a)\,a(b)\,\overline{a(c)\,a(d)}\,(abcd)^{-3/4}\,e\!\left(\frac{-(a+b-c-d)\,\overline h}k\right)\\
\qquad\cdot e\!\left(\mp\frac18\mp\frac18\mp\frac18\mp\frac18\right)\int\limits_M^{2M}x\,e\!\left(\pm\frac{2\sqrt{\vphantom bax}}k
\pm\frac{2\sqrt{bx}}k\pm\frac{2\sqrt{\vphantom bcx}}k\pm\frac{2\sqrt{dx}}k\right)\mathrm dx,
\end{multline*}
where $a$, $b$, $c$ and $d$ are positive integers from the interval $\left[1,N\right]$,
where the signs $\mp$ correspond to the $\pm$, and where otherwise all possible choices of the $\pm$-signs appear.

The main term of the fourth moment comes from the terms with two plus signs and two minus signs and in which $\sqrt{\vphantom ba}+\sqrt b=\sqrt{\vphantom bc}+\sqrt d$, or $\sqrt{\vphantom ba}+\sqrt{\vphantom bc}=\sqrt b+\sqrt d$ or $\sqrt{\vphantom ba}+\sqrt d=\sqrt b+\sqrt{\vphantom bc}$, depending on the locations of the minus signs. In this case the integral $\int_M^{2M}$ reduces to $3M^2/2$. The sum can be extended from $\sum_{a,b,c,d\leqslant N}$ to $\sum_{a,b,c,d=1}^\infty$ with an error
\[\ll k^2\,M^{2+\varepsilon}\,N^{-1/4},\]
as in (3.6) in \cite{Tsang1992}.

The terms with four plus signs or four minus signs are estimated by absolute values by Lemma \ref{fourth-moment-lemma} giving
\begin{align*}
&\ll k^2\sum_{a\leqslant N}\sum_{b\leqslant N}\sum_{c\leqslant N}\sum_{d\leqslant N}\left|a(a)\,a(b)\,a(c)\,a(d)\right|(abcd)^{-3/4}\cdot\frac{k\,M^{3/2}}{\sqrt{\vphantom ba}+\sqrt b+\sqrt{\vphantom bc}+\sqrt d}\\
&\ll k^3\,M^{3/2}\left(\sum_{n\leqslant N}\left|a(n)\right|\,n^{-7/8}\right)^{\!4}
\ll k^3\,M^{3/2}\,(N^{1/8})^4=k^3\,M^{3/2}\,N^{1/2}.
\end{align*}

Let us next consider the terms involving three signs of one kind and one of the other kind. Without loss of generality we may consider the signs $+$, $+$, $+$ and $-$. We first observe that the contribution of those terms in which $\Delta=\sqrt{\vphantom ba}+\sqrt b+\sqrt{\vphantom bc}-\sqrt d$ vanishes cancel simply because for them
\[e\!\left(+\frac18+\frac18+\frac18-\frac18\right)e\!\left(\frac{2\Delta\,\sqrt x}k\right)+
e\!\left(-\frac18-\frac18-\frac18+\frac18\right)e\!\left(\frac{-2\Delta\,\sqrt x}k\right)=0.\]

Next we split the sums $\sum_a\sum_b\sum_c\sum_d$ dyadically in each variable into $\ll\log^4N\ll M^\varepsilon$ subsums of the form
\[\sum_{A<a\leqslant2A}\sum_{B<b\leqslant2B}\sum_{C<c\leqslant2C}\sum_{D<d\leqslant2D},\]
where naturally $A$, $B$, $C$ and $D$ are all $\ll N$. It turns out, that by symmetry, it is enough to consider the terms in which $A\ll C$ and $B\ll C$. In particular, we do these simplifications in order to split the subsums further in terms of a dyadic decomposition of the range of $\Delta$ and apply Lemma \ref{counting-three-one-quadruples} to count the number of terms in each subsubsum.

We define $\delta$ to denote $C^{-1/2}\left|\Delta\right|$. The rest of estimating the $($$+$$+$$+$$-$$)$-terms is divided into three cases according to whether $\delta\gg1$, $1/C\ll\delta\ll1$ or $\delta\ll1/C$. The terms with $\delta\gg1$ are the easiest to dispose of. Namely, there are trivially at most $\ll ABCD$ such terms, and so by Lemma \ref{fourth-moment-lemma}, they contribute
\begin{align*}
&\ll k^2\,ABCD\,(ABCD)^{-3/4}\cdot\frac{k\,M^{3/2+\varepsilon}}{\delta\,\sqrt C}\\
&\ll k^3\,M^{3/2+\varepsilon}\,(ABCD)^{1/4}\,C^{-1/2}\\
&\ll k^3\,M^{3/2+\varepsilon}\,(CD)^{1/4}
\ll k^3\,M^{3/2+\varepsilon}\,N^{1/2}.
\end{align*}

Next, let us consider the terms with $1/C\ll\delta\ll1$. Since
\[\delta^2\,C=\bigl(\sqrt{\vphantom ba}+\sqrt b+\sqrt{\vphantom bc}-\sqrt d\,\bigr)^2
=\bigl(\sqrt{\vphantom ba}+\sqrt b-\sqrt d\,\bigr)^2-c+2\bigl(\sqrt{\vphantom ba}+\sqrt b+\sqrt{\vphantom bc}-\sqrt d\,\bigr)\sqrt{\vphantom bc},\]
we have
\[c=\bigl(\sqrt{\vphantom ba}+\sqrt b-\sqrt d\,\bigr)^2+O(\delta\,C),\]
and so there are $\ll\delta C$ possible values of $c$ for any given triple $\left\langle a,b,d\right\rangle$. Thus, there are $\ll AB\delta CD$ terms with $1/C\ll\delta\ll1$, and by Lemma \ref{fourth-moment-lemma}, they contribute
\[\ll k^2\,AB\delta CD\,(ABCD)^{-3/4}\cdot\frac{k\,M^{3/2+\varepsilon}}{\delta\,\sqrt C}\ll k^3\,M^{3/2+\varepsilon}\,N^{1/2}.\]

We next consider the third case in which $\delta\ll1/C$. This case will be split into two subcases depending on whether $C\gg\Gamma$ or $C\ll\Gamma$, where $\Gamma$ is some positive real number whose value will be set later when its impact on the final error terms is easier to see.

Let first $C\gg\Gamma$. We split our terms into further subsums by performing a dyadic division of the value range of $\delta$. Since by Lemma \ref{lower-bound-for-square-root-expression} we have
\[\delta=C^{-1/2}\left|\Delta\right|\gg C^{-5/2}\,(ABC)^{-1/2}\gg N^{-4},\]
there will be $\ll\log N\ll M^\varepsilon$ such subranges of $\delta$ to consider. Let us consider the terms where $\delta$ lies in one of these. Let us observe that $D\ll C$. By the first upper bound given by Lemma \ref{counting-three-one-quadruples} there are $\ll ABCD\,(\delta+C^{-3/2+\varepsilon})$ corresponding terms, which then contribute, estimating the integral $\int_M^{2M}$ either by absolute values or by Lemma \ref{fourth-moment-lemma},
\begin{align*}
&\ll k^2\,M^{\varepsilon}\,ABCD\left(\delta+C^{-3/2+\varepsilon}\right)
\left(ABCD\right)^{-3/4}\min\left\{M^2,\frac{k\,M^{3/2}}{\delta\,\sqrt C}\right\}\\
&\ll k^3\,M^{3/2+\varepsilon}\,(ABCD)^{1/4}\,C^{-1/2}
+k^2\,M^{2+\varepsilon}\,(ABCD)^{1/4}\,C^{-3/2+\varepsilon}\\
&\ll k^3\,M^{3/2+\varepsilon}\,N^{1/2}
+k^2\,M^{2+\varepsilon}\,\Gamma^{-1/2}.
\end{align*}

Next, let $C\ll\Gamma$. We again perform the same dyadic division of the value range of $\delta$. In one such subrange, the second upper bound of Lemma \ref{counting-three-one-quadruples} says that there $\ll ABCD\left(\delta C^2+(ABCD)^{-1/2}\right)C^\varepsilon$ corresponding terms which then contribute, estimating the integrals by Lemma \ref{fourth-moment-lemma} and $1/\delta$ by Lemma \ref{lower-bound-for-square-root-expression},
\begin{align*}
&\ll k^2\,M^\varepsilon\,ABCD\left(\delta C^2+(ABCD)^{-1/2}\right)C^\varepsilon
(ABCD)^{-3/4}\,\frac{k\,M^{3/2+\varepsilon}}{\delta\,\sqrt C}\\
&\ll k^3\,M^{3/2+\varepsilon}\left(C^{5/2}+C^{3/2}(ABCD)^{1/4}\right)
\ll k^3\,M^{3/2+\varepsilon}\,\Gamma^{5/2}.
\end{align*}
This concludes our treatment of the $($$+$$+$$+$$-$$)$-terms.

The last group of signs to consider are those with two plus signs and two minus signs. Those terms in which $\Delta=\sqrt{\vphantom ba}+\sqrt b-\sqrt{\vphantom bc}-\sqrt d$ vanishes were already considered in the derivation of the main term, and so we may assume throughout that $\Delta\neq0$. We perform again the dyadic division of the summations in $a$, $b$, $c$ and $d$, and by symmetry, we may focus on the terms with $A\ll C$, $B\ll C$ and $D\ll C$.

We may again define $\delta=C^{-1/2}\left|\Delta\right|$, and the regions $\delta\gg1$ and $1/C\ll\delta\ll1$ are handled by the same estimates as in the $($$+$$+$$+$$-$$)$-terms. The case $\delta\ll1/C$ is also similar, and split into subcases depending on whether $C\gg\Gamma$ or $C\ll\Gamma$. The latter works out in exactly the same as for the $($$+$$+$$+$$-$$)$-terms, except that Lemma \ref{counting-two-two-quadruples} is to be used instead of Lemma \ref{counting-three-one-quadruples}.

The only new complication arising from the case $C\gg\Gamma$ is that the first bound of Lemma \ref{counting-two-two-quadruples} has an extra term which does not appear in Lemma \ref{counting-three-one-quadruples}.
Since $\Delta\ll C^{-1/2}$, at least one of $A$ and $B$ must be $\asymp C$. Let us suppose that $B\asymp C$.
The contribution from the extra term is, estimating the integrals by absolute values and $AD\gg1$,
\begin{align*}
&\ll k^2\,M^\varepsilon\left(ABCD\right)^{1/2}\left(ABCD\right)^{-3/4}\,M^2\\
&\ll k^2\,M^{2+\varepsilon}\left(ABCD\right)^{-1/4}
\ll k^2\,M^{2+\varepsilon}\,C^{-1/2}
\ll k^2\,M^{2+\varepsilon}\,\Gamma^{-1/2}.
\end{align*}

Finally, we only need to collect all the error terms and choose suitable values for $N$ and $\Gamma$. We have established that
\begin{multline*}\int\limits_M^{2M}\left|A\!\left(x,\frac hk\right)\right|^4\mathrm dx
-C_F\,k^2\,M^2
\ll k^{5/2}\,M^{9/4+\varepsilon}\,N^{-1/2}
+k^4\,M^{3+\varepsilon}\,N^{-2}\\
+k^2\,M^{2+\varepsilon}\,N^{-1/4}
+k^3\,M^{3/2+\varepsilon}\,N^{1/2}
+k^3\,M^{3/2+\varepsilon}\,\Gamma^{5/2}
+k^2\,M^{2+\varepsilon}\,\Gamma^{-1/2}.\end{multline*}
Choosing $\Gamma=k^{-1/3}\,M^{1/6}$ optimizes the last two terms to $\ll k^{13/6}\,M^{23/12+\varepsilon}$. The value of $N$ can be chosen by optimizing the first and fourth terms; we choose $N=k^{-1/2}\,M^{3/4}$ and the first and fourth terms simplify to $\ll k^{11/4}\,M^{15/8+\varepsilon}$. The second and third terms will become $k^5\,M^{3/2+\varepsilon}+k^{9/8}\,M^{29/16+\varepsilon}$. This is easily seen to be $\ll k^{13/6}\,M^{23/12+\varepsilon}+k^{11/4}\,M^{15/8+\varepsilon}$, and we are done.

\section{More moment estimates}

In the proof of Theorem \ref{oscillation-result}, we will need a lower bound for the mean square of $A_{++}$. We will prove this following Ivi\'c and Zhai \cite{IvicEtZhai2012}.

\begin{proposition}\label{primitive}
Let $M\geqslant1$, let $0\leqslant\Xi\leqslant M$, and let $h$ and $k$ be coprime integers with $1\leqslant k\leqslant M$. Then
\[\int\limits_M^{M+\Xi}A\!\left(x,\frac hk\right)\mathrm dx
\ll k^{3/2}\,M^{3/4}+k\,\Xi\,M^\varepsilon.\]
\end{proposition}

Using the truncated Voronoi identity, i.e.\ Theorem \ref{truncated-voronoi-identity} above, and Lemma~\ref{primitive-lemma} we obtain
\begin{align*}
&\int\limits_M^{M+\Xi}A\!\left(x,\frac hk\right)\mathrm dx\\
&=\frac{k^{1/2}}{\pi\,\sqrt2}
\int\limits_M^{M+\Xi}x^{1/4}\sum_{n\leqslant M}a(n)\,n^{-3/4}\,e\!\left(\frac{-n\overline h}k\right)\cos\!\left(\frac{4\pi\sqrt{nx}}k-\frac\pi4\right)\mathrm dx\\
&\qquad+\int\limits_M^{M+\Xi}O\!\left(k\,M^\varepsilon\right)\mathrm dx\\
&=\frac{k^{1/2}}{\pi\,\sqrt2}\sum_{n\leqslant M}a(n)\,n^{-3/4}\,e\!\left(\frac{-n\overline h}k\right)\int\limits_M^{M+\Xi}x^{1/4}\,\cos\!\left(\frac{4\pi\sqrt{nx}}k-\frac\pi4\right)\mathrm dx\\
&\qquad+O\!\left(k\,\Xi\,M^\varepsilon\right)\\
&\ll k^{1/2}\sum_{n\leqslant M}\left|a(n)\right|n^{-3/4}\cdot\frac k{\sqrt n}\,M^{3/4}+k\,\Xi\,M^\varepsilon
\ll k^{3/2}\,M^{3/4}+k\,\Xi\,M^\varepsilon.
\end{align*}

\begin{proposition}\label{first-moment}
Let $M\gg1$, and let $h$ and $k$ be coprime integers with $1\leqslant k\ll M^{1/6-\varepsilon}$. Then
\[\int\limits_M^{2M}\left|A\!\left(x,\frac hk\right)\right|\mathrm dx
\gg k^{1/2}\,M^{5/4}.\]
\end{proposition}

Using Theorem \ref{mean-square-of-long-sums}, the Cauchy--Schwarz inequality, and Corollary \ref{cubic-moment} (remembering that Theorem \ref{fourth-moment-estimate} removes the $\varepsilon$) we get
\begin{align*}
k\,M^{3/2}&\asymp\int\limits_M^{2M}\left|A\!\left(x,\frac hk\right)\right|^2\mathrm dx\\
&\ll\sqrt{\int\limits_M^{2M}\left|A\!\left(x,\frac hk\right)\right|^3\mathrm dx}
\,\sqrt{\int\limits_M^{2M}\left|A\!\left(x,\frac hk\right)\right|\mathrm dx}\\
&\ll\sqrt{k^{3/2}\,M^{7/4}}\,\sqrt{\int\limits_M^{2M}\left|A\!\left(x,\frac hk\right)\right|\mathrm dx},
\end{align*}
and the claim follows easily.

\bigbreak
\begin{lemma}\label{first-quadrant-lower-bound}
Let $M\gg1$, and let $h$ and $k$ be coprime integers with $1\leqslant k\ll M^{1/6-\varepsilon}$. Then
\[\int\limits_M^{2M}\left|A_{++}\!\left(x,\frac hk\right)\right|^2\mathrm dx
\gg k\,M^{3/2}.\]
\end{lemma}

By Proposition \ref{first-moment} and Proposition \ref{primitive}, we have
\begin{align*}
k^{1/2}\,M^{5/4}
&\ll\int\limits_M^{2M}\left|A\!\left(x,\frac hk\right)\right|\mathrm dx\\
&\ll\int\limits_M^{2M}\left|\Re A\!\left(x,\frac hk\right)\right|\mathrm dx
+\int\limits_M^{2M}\left|\Im A\!\left(x,\frac hk\right)\right|\mathrm dx\\
&=2\int\limits_M^{2M}\left(\left(\Re A\!\left(x,\frac hk\right)\!\right)_{\!+}
+\left(\Im A\!\left(x,\frac hk\right)\!\right)_{\!+}\right)\mathrm dx\\
&\qquad-\int\limits_M^{2M}\left(\Re A\!\left(x,\frac hk\right)
+\Im A\!\left(x,\frac hk\right)\!\right)\mathrm dx
\end{align*}
For sufficiently large $M$, the last integral, which is $\ll k^{3/2}\,M^{3/4}+k\,M^{1+\varepsilon}$, may be absorbed to the left-hand side, and we may continue the argument with the Cauchy--Schwarz inequality:
\begin{align*}
k^{1/2}\,M^{5/4}
&\ll\int\limits_M^{2M}\left(\left(\Re A\!\left(x,\frac hk\right)\!\right)_{\!+}
+\left(\Im A\!\left(x,\frac hk\right)\!\right)_{\!+}\right)\mathrm dx\\
&\ll\int\limits_M^{2M}\left|A_{++}\!\left(x,\frac hk\right)\right|\mathrm dx
\ll\sqrt M\,\sqrt{\int\limits_M^{2M}\left|A_{++}\!\left(x,\frac hk\right)\right|^2\mathrm dx},
\end{align*}
and the result follows easily.

\section{Proof of Theorem \ref{mean-square-of-very-short-sums}}

Choosing $N=M\asymp x$ in the truncated Voronoi identity, i.e.\ Theorem \ref{truncated-voronoi-identity} above, gives
\begin{align*}
&\sum_{x\leqslant n\leqslant x+\Delta}a(n)\,e\!\left(\frac{nh}k\right)\\
&\qquad=\frac{k^{1/2}}{\pi\sqrt2}\sum_{n\leqslant M}a(n)\,n^{-3/4}\,e\!\left(\frac{-n\overline h}k\right)\\
&\qquad\qquad\cdot\Bigg(\cos\left(\frac{4\pi\sqrt{n\left(x+\Delta\right)}}k-\frac\pi4\right)\left(x+\Delta\right)^{1/4}\\
&\hspace*{14em}-\cos\left(\frac{4\pi\sqrt{nx}}k-\frac\pi4\right)x^{1/4}\Bigg)+O(kx^\varepsilon)\\
&\qquad=\frac{k^{1/2}\,x^{1/4}}{\pi\sqrt2}
\sum_{n\leqslant M}a(n)\,n^{-3/4}\,e\!\left(\frac{-n\overline h}k\right)\\
&\qquad\qquad
\cdot\left(\cos\left(\frac{4\pi\sqrt{n\left(x+\Delta\right)}}k-\frac\pi4\right)
-\cos\left(\frac{4\pi\sqrt{nx}}k-\frac\pi4\right)\right)+O(kM^\varepsilon).
\end{align*}
Here the last estimate follows from estimates by absolute values and the simple observation that
\[\left(x+\Delta\right)^{1/4}-x^{1/4}
=\frac14\int\limits_x^{x+\Delta}t^{-3/4}\,\mathrm dt
\asymp\Delta\,M^{-3/4}.\]
Plugging this output of the truncated Voronoi formula into the mean square expression gives
\begin{align*}
&\int\limits_M^{M+\Xi}\,\left|\sum_{x\leqslant n\leqslant x+\Delta}a(n)\,e\!\left(\frac{nh}k\right)\right|^2\mathrm dx\\
&=\frac{k}{2\pi^2}\int\limits_M^{M+\Xi}x^{1/2}
\left|\sum_{n\leqslant M}a(n)\,n^{-3/4}\,e\!\left(\frac{-n\overline h}k\right)\,\left(\cos\left(\ldots\right)-\cos\left(\cdots\right)\right)\right|^2\mathrm dx\\
&
+k^{1/2}\,\Re\int\limits_M^{M+\Xi}x^{1/4}
\sum_{n\leqslant M}a(n)\,n^{-3/4}\,e\!\left(\frac{-n\overline h}k\right)
\left(\cos\left(\ldots\right)-\cos\left(\ldots\right)\right)O\!\left(k\,M^\varepsilon\right)\mathrm dx\\
&\qquad\qquad\qquad+O\!\left(k^2\,\Xi\,M^\varepsilon\right).
\end{align*}
The last term is $\ll\Xi\,\Delta\,M^{-\varepsilon}$, and
once we have proved that the first term is $\ll\Xi\,\Delta$, an application of the Cauchy--Schwarz inequality gives the conclusion that the second term is
\[\ll\sqrt{\Xi\Delta}\sqrt{\Xi\,k^2\,M^\varepsilon}
=\Xi\,k\,\Delta^{1/2}\,M^{\varepsilon}
\ll\Xi\,\Delta\,M^{-\varepsilon}.\]

The rest of the proof consists basically of expanding the square, applying the first derivative test and collecting terms.
The main contribution comes from the low-frequency diagonal terms, which are easy enough to handle. However, one must be careful in order to avoid resonances in the off-diagonal terms.

First, let us introduce some notation: for $\lambda\in\mathbb R$, we shall write
\[s_\lambda
=\sum_{n\leqslant M}a(n)\,n^{-3/4}\,\cos\frac{2\pi n\overline h}k
\,e\!\left(\frac{2\sqrt{n\left(x+\lambda\right)}}k-\frac18\right),\]
where $\lambda$ will ultimately be either $0$ or $\Delta$.
Similarly, we shall write $\widetilde s_\lambda$ for the same sum with $\cos\left(2\pi n\overline h/k\right)$ replaced by $\sin\left(2\pi n\overline h/k\right)$. Thus we shall have
\begin{align*}
&\int\limits_M^{M+\Xi}\,\left|\sum_{x\leqslant n\leqslant x+\Delta}a(n)\,e\!\left(\frac{nh}k\right)\right|^2\mathrm dx\\
&\qquad=\frac k{2\pi^2}\int\limits_M^{M+\Xi}x^{1/2}\left(\Re s_\Delta-\Re s_0\right)^2\mathrm dx\\
&\qquad\qquad+\frac k{2\pi^2}\int\limits_M^{M+\Xi}x^{1/2}\left(\Re\widetilde s_\Delta-\Re\widetilde s_0\right)^2\mathrm dx
+O\!\left(\Xi\,\Delta\,M^{-\varepsilon}\right).
\end{align*}
The second integral on the right-hand side can be treated in exactly the same way as the first one. In the end the factors $\cos\left(2\pi n\overline h/k\right)$ and $\sin\left(2\pi n\overline h/k\right)$ will go away due to the fact that $\cos^2\alpha+\sin^2\alpha=1$. In the following we shall, for the sake of simpler notation, to assume that the coefficients $a(n)$ are real. If this is not the case, then we can just replace $a(n)\cos2\pi n\overline h/k$ and $a(n)\sin2\pi n\overline h/k$ by the real and imaginary parts of $a(n)\,e(nh/k)$.

The integrand $\left(\Re s_\Delta-\Re s_0\right)^2$ will be exchanged for a handful of other terms which lead to exponential integrals which can be fed to the first derivative test. The first step will be splitting the sums into low-frequency and high-frequency terms:
\[s_\lambda^{\leqslant}=\sum_{n\leqslant N_0}\ldots,\qquad
s_\lambda^>=\sum_{N_0<n\leqslant M}\ldots,\]
where $N_0=M/(2\Delta)$. Furthermore, we shall denote the corresponding real parts by $\sigma_\lambda^\leqslant$ and $\sigma_\lambda^>$, and the imaginary parts by $t_\lambda^\leqslant$ and $t_\lambda^>$.

The low- versus high-frequency split will then yield
\begin{multline*}
\left(\Re s_\Delta-\Re s_0\right)^2
=\left((\sigma_\Delta^\leqslant-\sigma_0^\leqslant)+(\sigma_\Delta^>-\sigma_0^>)\right)^2\\
=(\sigma_\Delta^\leqslant-\sigma_0^\leqslant)^2
+(\sigma_\Delta^>-\sigma_0^>)^2
+2\,(\sigma_\Delta^\leqslant-\sigma_0^\leqslant)\,(\sigma_\Delta^>-\sigma_0^>).
\end{multline*}
These three terms and their integrals will be handled separately. We start with the last one.

The last term can be written in terms of the exponential function:
\begin{multline*}
2\,(\sigma_\Delta^\leqslant-\sigma_0^\leqslant)\,(\sigma_\Delta^>-\sigma_0^>)
\ll2\,(s_\Delta^\leqslant-s_0^\leqslant)\,(\sigma_\Delta^>-\sigma_0^>)\\
=(s_\Delta^\leqslant-s_0^\leqslant)\,\overline{(s_\Delta^>-s_0^>)}+(s_\Delta^\leqslant-s_0^\leqslant)\,(s_\Delta^>-s_0^>).
\end{multline*}
By Lemmas \ref{mean-square-lemma-2} and \ref{mean-square-lemma-3}, the contribution from the integral of the first term will be
\begin{align*}
&\ll k\int\limits_M^{M+\Xi}x^{1/2}\,(s_\Delta^\leqslant-s_0^\leqslant)\,\overline{(s_\Delta^>-s_0^>)}\,\mathrm dx\\
&\ll k\sum_{m\leqslant N_0}\sum_{N_0<n\leqslant M}
\left|a(m)\,a(n)\right|\left(mn\right)^{-3/4}\cdot\frac{k\,\sqrt{n}\,M}{n-m}\\
&\ll k^2\,M^{1+\varepsilon}
\ll \Xi\,\Delta\,M^{-\varepsilon}.
\end{align*}
The contribution from the other term, the one without the complex conjugation, can be estimated in the same way using Lemma~\ref{mean-square-lemma-1}.

The middle term can be estimated by two mean squares:
\[(\sigma_\Delta^>-\sigma_0^>)^2
=(\sigma_\Delta^>)^2+(\sigma_0^>)^2-2\,\sigma_\Delta^>\,\sigma_0^>
\ll(\sigma_\Delta^>)^2+(\sigma_0^>)^2
\ll\bigl|s_\Delta^>\bigr|^2+\bigl|s_0^>\bigr|^2.
\]
The contribution from the integral of either of these is obtained by expanding the square, considering the diagonal and the non-diagonal terms separately, estimating the diagonal terms by absolute values, and estimating the non-diagonal terms using Lemma \ref{mean-square-lemma-2}, giving
\begin{align*}
\frac k{2\pi^2}\int\limits_M^{M+\Xi}x^{1/2}\,\bigl|s_\lambda^>\bigr|^2\,\mathrm dx
&\ll k\sum_{N_0<n\leqslant M}\left|a(n)\right|^2n^{-3/2}\,\Xi\,M^{1/2}\\
&\quad+k\sum_{N_0<m\leqslant M}\sum_{m<n\leqslant M}\left|a(m)\,a(n)\right|(mn)^{-3/4}\cdot\frac{k\,\sqrt{n}\,M}{n-m}\\
&\ll k\,\Xi\,M^{1/2}\,N_0^{-1/2}+k^2\,M^{1+\varepsilon}\\
&\ll k\,\Xi\,\Delta^{1/2}+\Xi\,\Delta\,M^{-\varepsilon}
\ll\Xi\,\Delta\,M^{-\varepsilon}.
\end{align*}

The first term will give the main contribution. First, we split it into two parts, the second of which will give rise only to oscillating integrals:
\begin{align*}
(\sigma_\Delta^\leqslant-\sigma_0^\leqslant)^2
&=\frac12\left((\sigma_\Delta^\leqslant-\sigma_0^\leqslant)^2+(t_\Delta^\leqslant-t_0^\leqslant)^2\right)+\frac12\left((\sigma_\Delta^\leqslant-\sigma_0^\leqslant)^2-(t_\Delta^\leqslant-t_0^\leqslant)^2\right)\\
&=\frac12\bigl|s_\Delta^\leqslant-s_0^\leqslant\bigr|^2
+\frac12\Re(s_\Delta^\leqslant-s_0^\leqslant)^2.
\end{align*}
The integral of the $O$-term can be estimated by expanding the square and estimating the resulting exponential integrals by Lemma \ref{mean-square-lemma-1}:
\begin{align*}
\frac k{2\pi^2}\int\limits_M^{M+\Xi}x^{1/2}\,(s_\Delta^\leqslant-s_0^\leqslant)^2\,\mathrm dx
&\ll k\sum_{m\leqslant N_0}\sum_{n\leqslant N_0}
\left|a(m)\,a(n)\right|(mn)^{-3/4}\cdot\frac{k\,M}{\sqrt n+\sqrt m}\\
&\ll k^2\,M^{1+\varepsilon}
\ll\Xi\,\Delta\,M^{-\varepsilon}.
\end{align*}

The last integral we have to consider is
\[\frac k{4\pi^2}\int\limits_M^{M+\Xi}x^{1/2}\,\bigl|s_\Delta^\leqslant-s_0^\leqslant\bigr|^2\,\mathrm dx.\]
Here we again expand the square and consider the diagonal and the non-diagonal terms separately.
The non-diagonal terms are estimated as before using Lemma \ref{mean-square-lemma-3} giving
\begin{align*}
\ll k\sum_{m\leqslant N_0}\sum_{m<n\leqslant N_0}
\left|a(m)\,a(n)\right|(mn)^{-3/4}\frac{k\,M\,\sqrt n}{n-m}
\ll k^2\,M^{1+\varepsilon}
\ll\Xi\,\Delta\,M^{-\varepsilon}.
\end{align*}

The diagonal terms give the main contribution; they are
\begin{align*}
\frac k{4\pi^2}\sum_{n\leqslant N_0}\frac{\left|a(n)\right|^2}{n^{3/2}}
\cos^2\frac{2\pi n\overline h}k
\int\limits_M^{M+\Xi}x^{1/2}
\left|e\!\left(\frac{2\sqrt n\left(\sqrt{x+\Delta}-\sqrt x\right)}k\right)-1\right|^2\mathrm dx.
\end{align*}
The integrals involving $\widetilde s_\lambda$ give the same main terms with $\cos^2$ replaced by $\sin^2$, and so these trigonometric factors cancel away leaving only
\begin{align*}
\frac k{4\pi^2}\sum_{n\leqslant N_0}\frac{\left|a(n)\right|^2}{n^{3/2}}
\int\limits_M^{M+\Xi}x^{1/2}
\left|e\!\left(\frac{2\sqrt n\left(\sqrt{x+\Delta}-\sqrt x\right)}k\right)-1\right|^2\mathrm dx.
\end{align*}
This sum is split into low-frequence terms and high-frequency terms according to whether $n\leqslant N_1$ or $N_1<n\leqslant N_0$, where $N_1=\frac14\,k^2\,M\,\Delta^{-2}$. The high-frequency terms give a nonnegative contribution which is at most
\[\ll k\sum_{N_1<n\leqslant N_0}\frac{\left|a(n)\right|^2}{n^{3/2}}\,\Xi\,M^{1/2}
\ll k\,N_1^{-1/2}\,\Xi\,M^{1/2}
\ll\Xi\,\Delta.\]
In the low-frequency terms with $n\leqslant N_1$, we have
\[0<\frac{2\sqrt n\left(\sqrt{x+\Delta}-\sqrt x\right)}k
=\frac{\sqrt n}k\int\limits_x^{x+\Delta}\frac{\mathrm dt}{\sqrt t}
\leqslant\frac{\sqrt n\,\Delta}{k\,\sqrt M}\leqslant\frac12,\]
where $x\in\left[M,M+\Xi\right]$. Thus, the low-frequency terms are
\begin{align*}
&\asymp k\sum_{n\leqslant N_1}\frac{\left|a(n)\right|^2}{n^{3/2}}
\,\Xi\,M^{1/2}\,\frac{n\,\Delta^2}{k^2\,M}
\ll k^{-1}\,\Xi\,M^{-1/2}\,\Delta^2\,N_1^{1/2}
\asymp\Xi\,\Delta,
\end{align*}
and we are done.

By inspecting the last line, we also observe that, if $\Delta\ll M^{1/2-\varepsilon}$ and $M\gg_\varepsilon1$, then $N_1\gg M^\varepsilon$ and the low-frequency terms actually are $\asymp\Xi\,\Delta$, so that we get the second conclusion of the theorem.

\section{Proof of Theorem \ref{sharper-mean-square-for-short-sums}}

Let $\lambda\in\mathbb Z_+$ be such that $2^\lambda\leqslant\Delta^{1/2}<2^{\lambda+1}$, and write $b=\Delta\,2^{-\lambda}$. Then $b\asymp 2^\lambda\asymp\Delta^{1/2}$. The idea of the proof is to consider subsums of length $b$ instead of individual terms. The relevant observation here is that
\begin{multline*}
\max_{1\leqslant U\leqslant\Delta}
\left|\sum_{x\leqslant n\leqslant x+U}a(n)\,e\!\left(\frac{nh}k\right)\right|^2\\
\ll\max_{1\leqslant j<2^{\lambda}}
\left|\sum_{x<n\leqslant x+jb}a(n)\,e\!\left(\frac{nh}k\right)\right|^2
+O\!\left(b^2\log^2M\right),
\end{multline*}
where $j$ takes only integral values and we have used Deligne's and Shiu's estimates (see Subsection \ref{results-on-holomorphic-cusp-forms} above).

The maximum is certainly attained for some value particular value of $j$ which we shall call $j_0$. The remainder of the argument is really just a matter of estimating the sum of length $j_0b$ so that the dependence on $j_0$ goes away, and then finishing off with Theorem \ref{mean-square-of-very-short-sums}.

We use the binary representation of $j_0$ to dyadically dissect $\left[x,x+j_0b\right]$. Write
\[j_0=2^{\lambda_1}+2^{\lambda_2}+\ldots+2^{\lambda_N},\]
where $N\in\mathbb Z_+$ and the exponents are integers satisfying
\[0\leqslant\lambda_N<\ldots<\lambda_2<\lambda_1<\lambda.\]
Writing also
\[\Lambda_0=0,\quad
\Lambda_1=2^{\lambda_1},\quad
\Lambda_2=2^{\lambda_1}+2^{\lambda_2},\quad
\ldots,\quad
\Lambda_N=2^{\lambda_1}+2^{\lambda_2}+\ldots+2^{\lambda_N},\]
we estimate
\[\left|\sum_{x<n\leqslant j_0b}a(n)\,e\!\left(\frac{nh}k\right)\right|^2
\ll N
\sum_{k=1}^N
\left|\sum_{x+\Lambda_{k-1}b<n\leqslant x+\Lambda_kb}
a(n)\,e\!\left(\frac{nh}k\right)\right|^2.\]

Writing next $\Lambda_k=\nu_k\,2^{\lambda_k}$ for $k\in\left\{0,1,\ldots,N\right\}$, where we set $\lambda_0=0$, we have $0\leqslant\nu_k<2^{\lambda-\lambda_k}$, and we may estimate
\begin{align*}
&\left|\sum_{x+\Lambda_{k-1}b<n\leqslant x+\Lambda_kb}a(n)\,e\!\left(\frac{nh}k\right)\right|^2\\
&\qquad=\left|\sum_{x+\nu_{k-1}2^{\lambda_{k-1}}b<n\leqslant x+\nu_{k-1}2^{\lambda_{k-1}}b+2^{\lambda_k}b}a(n)\,e\!\left(\frac{nh}k\right)\right|^2\\
&\qquad\leqslant
\sum_{\nu=0}^{2^{\lambda-\lambda_{k-1}}-1}
\left|\sum_{x+\nu2^{\lambda_{k-1}}b<n\leqslant x+\nu2^{\lambda_{k-1}}b+2^{\lambda_k}b}a(n)\,e\!\left(\frac{nh}k\right)\right|^2.
\end{align*}
Combining this with $\int_M^{2M}$ and applying Theorem \ref{mean-square-of-very-short-sums} gives
\begin{align*}
&\int\limits_M^{2M}\max_{0\leqslant U\leqslant\Delta}\left|\sum_{x\leqslant n\leqslant x+U}a(n)\,e\!\left(\frac{nh}k\right)\right|^2\mathrm dx\\
&\qquad\ll N\sum_{k=1}^N
\sum_{\nu=0}^{2^{\lambda-\lambda_{k-1}}-1}
\int\limits_M^{2M}\,
\left|\sum_{x+\nu2^{\lambda_{k-1}}b<n\leqslant x+\nu2^{\lambda_{k-1}}b+2^{\lambda_k}b}a(n)\,e\!\left(\frac{nh}k\right)\right|^2\mathrm dx
\\
&\qquad\qquad+M\,\Delta\,\log^2M,
\end{align*}
and we may finish by estimating the first term on the right-hand side
\begin{align*}
&\qquad\ll\lambda\sum_{\ell=0}^{\lambda-1}
\sum_{\nu=0}^{2^{\lambda-(\ell+1)}-1}
\int\limits_M^{2M}\,
\left|\sum_{x+\nu2^{\ell+1}b<n\leqslant x+\nu2^{\ell+1}b+2^{\ell}b}a(n)\,e\!\left(\frac{nh}k\right)\right|^2\mathrm dx\\
&\qquad\ll\lambda\sum_{\ell=0}^{\lambda-1}
\sum_{\nu=0}^{2^{\lambda-(\ell+1)}-1}
M\,2^\ell\,b
\ll\lambda\sum_{\ell=1}^\lambda
2^{\lambda-\ell}\,M\,2^{\ell}\,b
\ll\lambda^2\,M\,2^\lambda\,b\ll M\,\Delta\,\log^2\Delta.
\end{align*}

\section{Proof of Theorem \ref{oscillation-result}}

We will only prove the case in which the argument of $A\!\left(x,h/k\right)$ is required to lie in the interval $\left[-\pi/2-\varepsilon/2,\pi+\varepsilon/2\right]$. The more general case follows from multiplying the underlying cusp form by a suitable unimodular constant.
We write $\Delta=c\,M^{1/2}\,\log^{-2}M$, and define a function
\[\omega(x)=\left|A_{++}\!\left(x,\frac hk\right)\right|^2-C\max_{0\leqslant U\leqslant\Delta}\left|A\!\left(x+U,\frac{h}k\right)-A\!\left(x,\frac{h}k\right)\right|^2-c\,k\,x^{1/2}\]
for $x\in\left[M,2M\right]$. Here $c$ is a small positive real constant and $C$ is a very large positive real constant. The point of this definition is that when $\omega(x)>0$ then $x$ must lie in a subinterval of length $\Delta$ in which $A\!\left(x,h/k\right)\gg k^{1/2}\,M^{1/4}$ and $A\!\left(x,h/k\right)$ has a positive real or positive imaginary part. Furthermore, since the value of $A(x+U,h/k)$, where $0\leqslant U\leqslant\Delta$, can deviate by at most $1/\sqrt C$ of the absolute value of $A(x,h/k)$ to the third quadrant of the complex plane, it is geometrically clear that if $C$ is sufficiently large depending on $\varepsilon$, say $C=1/\sin^2(\varepsilon/2)$, then the argument of $A(x+U,h/k)$ lies in the interval $\left[-\pi/2-\varepsilon/2,\pi+\varepsilon/2\right]$.

If $c$ is small enough, then we may apply Lemma \ref{first-quadrant-lower-bound} and Theorem~\ref{sharper-mean-square-for-short-sums} to get
\begin{align*}
\int\limits_M^{2M}\omega(x)\,\mathrm dx
\gg k\,M^{3/2}-O(C\,M\,\Delta\,\log^2M)-O(c\,k\,M^{3/2})
\gg k\,M^{3/2},
\end{align*}
and so $\int_M^{2M}\omega(x)\,\mathrm dx$ will be positive and $\gg k\,M^{3/2}$ for sufficiently large $M$.

Finally, let $S$ be the set of $x\in\left[M,2M\right]$ for which $\omega(x)>0$, and write $\ell$ for the measure of $S$. Then, applying Theorem \ref{fourth-moment-estimate},
\begin{align*}
k\,M^{3/2}
&\ll\int\limits_M^{2M}\omega(x)\,\mathrm dx
\ll\int\limits_S\omega(x)\,\mathrm dx
\ll\int\limits_S\left|A_{++}\!\left(x,\frac{h}k\right)\right|^2\mathrm dx\\
&\ll\sqrt{\int\limits_S1^2\,\mathrm dx}\,
\sqrt{\int\limits_M^{2M}\left|A\!\left(x,\frac{h}k\right)\right|^4\mathrm dx}
\ll\ell^{1/2}\,k\,M,
\end{align*}
so that $\ell\gg M$, and we are done.

\section*{Acknowledgements}

The author would like to express his gratitude for 
the valuable advice of A.-M. Ernvall-Hyt\"onen, for an insightful conversation with K.-M. Tsang on matters related to Voronoi-type formulae, and for the beneficial suggestions of an anonymous referee, especially regarding Theorem~\ref{large-values}.

The author first learned about the sign change results for the divisor problem and related problems at the excellent conference Elementare und analytische Zahlentheorie, held in Schloss Schney, August 13--18, 2012.
The author is grateful for the generous support of the organizers.

This research was funded by Finland's Ministry of Education through the Doctoral Program in Inverse Problems, the Academy of Finland through the Finnish Centre of Excellence in Inverse Problems Research, and the Foundation of Vilho, Yrj\"o and Kalle V\"ais\"al\"a.

\end{document}